\documentclass[12pt, oneside,reqno]{amsproc}   	
\usepackage{geometry}                		
\usepackage{graphicx}				
\usepackage{amssymb}
\usepackage{amsmath}
\usepackage[all]{xy}
\usepackage{mathtools}
\usepackage{tikz-cd}
\usepackage{enumitem}
\usepackage{tikz-cd}
\usepackage{setspace}

\newtheorem{thm}{Theorem}[section]
\newtheorem{lem}[thm]{Lemma}
\newtheorem{cor}[thm]{Corollary}
\newtheorem{prop}[thm]{Proposition}

\newtheorem{defn}{Definition}[section]

\numberwithin{equation}{section}

\renewcommand{\l}{\lambda}

\def\Pb{\ifmmode{\Bbb P}\else{$\Bbb P$}\fi}
\def\Z{\ifmmode{\Bbb Z}\else{$\Bbb Z$}\fi}
\def\C{\ifmmode{\Bbb C}\else{$\Bbb C$}\fi}
\def\R{\ifmmode{\Bbb R}\else{$\Bbb R$}\fi}
\def\S{\ifmmode{S^2}\else{$S^2$}\fi}

\def\S{\cal S}

\newenvironment{pf}{\paragraph{Proof:}}{\hfill$\square$ \newline}

\begin {document}
	
\title{A finiteness theorem via the mean curvature flow with surgery}
\begin{abstract}In this article, we use the recently developed mean curvature flow with surgery for 2-convex hypersurfaces to prove certain isotopy existence and finally extrinsic finiteness results (in the spirit of Cheeger's finiteness theorem) for the space of 2-convex closed embedded hypersurfaces in $\R^{n+1}$.  \end{abstract}

\author {Alexander Mramor}
\address{Department of Mathematics, University of California Irvine, CA 92617}
\email{mramora@uci.edu}

\maketitle

\section{Introduction}
Geometric flows, roughly speaking, allow us to deform the geometry of certain classes of manifolds (depend on how the flow is defined) and thus provides a tool to understand that class of manifolds. For example, a classical result of Huisken \cite{Hui} tells us that under the mean curvature flow all convex hypersurfaces (that is, hypersurfaces where all the principal curavtures at every point are positive) flow to round spheres - phrased another way this tells us that the set of all convex hypersurfaces (of $\R^{n+1}$) is connected, with the flow providing a canonical path of sorts between any two of the manifolds. By the smooth dependence of the flow on initial conditions, it is even a retraction. 
$\medskip$

However it's not so easy to get the maximal amount of information possible using a particular flow (or even if much information). A conjecture of Smale (later proved by Hatcher using non flow arguments in \cite{Hat}) tells us roughly that any continuous family of embedded $S^2 \subset \R^3$ (in smooth topology) can be deformed continuously to round $S^2$. An example of such an embedded sphere is a dumbbell, two copies of round $S^2$ connected by a thin neck; if the neck is sufficiently thin then under the mean curvature flow it shrinks much faster than the spheres and the flow can't be continued past some time $T_0$ because it develops a singularity:
 \begin{center}
$\includegraphics[scale = .3]{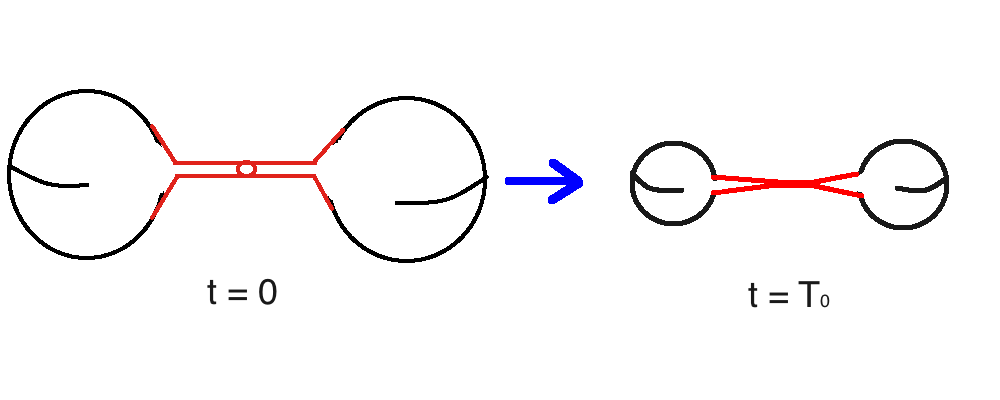}$
\end{center}
So the flow doesn't provide an isotopy of the dumbbell to a round sphere, even though there clearly is one. 
$\medskip$

As a way to continue the flow past singularities like this, first Huisken, Sinestari, and Brendle (specifically Brendle and Huisken tackled the $n=2$ case) and later independently Haslhofer and Kleiner developed a mean curvature flow with surgery for 2-convex embedded $n$-hypersurfaces $\mathcal{M}$ of $\R^{n+1}$. When applying the mean curvature flow on 2-convex hypersurfaces, it turns out (speaking for the Huisken-Sinestrari-Brendle approach) that the high curvature regions are either contained in necks possibly with caps or the hypersurface is uniformly convex. Doing surgery on the necks where the curvature is relatively low will then separate lower curvature regions from higher curvature regions which we can identify either as copies of $S^n$ or $S^{n-1} \times S^1$. Then we can continue the flow on the low curvature regions and repeat the process. Whittling away at the hypersurface this way until nothing is left the mentioned authors showed: 

\begin{thm}(Corollary 1.2 in \cite{HS2} for $n \geq 3$, but by \cite{BH} also true for $n=2$) Any smooth $n$-dimensional 2-convex immersed surface $F_0: \mathcal{M} \to \R^{n+1}$ with $n \geq 2$ is diffeomorphic to a finite connect sum of $S^{n-1} \times S^1$. Furthermore there exists a handlebody $\Omega$ and an immersion $G: \Omega \to \R^{n+1}$ such that $\partial \Omega \cong \mathcal{M}$ is diffeomorphic to the initial hypersurface $\mathcal{M}$ and such that $G|_{\partial \Omega} = F_0$. 
\end{thm}
In particular, one obtains a Schoenflies type theorem for simply connected 2-convex surfaces:
\begin{thm}(Corollary 1.3 in \cite{HS2} for $n \geq 3$, but by \cite{BH} also true for $n=2$) Any smooth closed simply connected $n$-dimensional 2-convex embedded surface $\mathcal{M} \subset \R^{n+1}$ with $n \geq 3$ is diffeomorphic to $S^n$ and bounds a region whose closure is diffeomorphic to a smoothly embedded $(n+1)$-dimensional standard closed ball. 
\end{thm}

These results of course are very nice but aren't quite in the same spirit as Huisken's or Hatcher's results; surgery is a discontinuous and paths of hypersurfaces aren't being created. But by essentially ``extending surgery necks by hand" analogous to \cite{Mar} we can use the flow with surgery to attain certain paths to canonical (or at least more canonical) representatives in the relevant isotopy classes. The cresecendo of the article is the finiteness theorem (corollary 1.6) below so to greatly simplify our arguments we focus on attaining monotonicity of paths in the theorems below instead of preserving 2-convexity, which has algready in large part been carried out in \cite{BHH}, \cite{BHH1} (see below for more details). Our results build one on after the other and are listed in logical order: 
\begin{thm} $\textbf{(Path to round sphere)}$ Any 2-convex hypersurface $\mathcal{M}^n \subset \R^{n+1}$ diffeomorphic to $S^n$ is isotopic to a round sphere through a monotone isotopy \end{thm}
By monotone isotopy we mean that if $H: M \times [0,1] \to \R^{n+1}$ is an (ambient) isotopy from $\mathcal{M}_0$ to $\mathcal{M}_1$ then it is monotone if (the hull of) $\mathcal{M}_s \subset \mathcal{M}_t$ for $t \leq s$ (also throughout isotopies are taken to be smooth in time unless otherwise indicated). Using this theorem we next show:
\begin{thm} $\textbf{(Torus to knot)}$ Given any 2-convex closed hypersurface $\mathcal{M}^n \subset \R^{n+1}$ diffeomorphic to $S^{1} \times S^{n-1}$ there exists a knot $\gamma$ so that for all $\epsilon > 0$ $\mathcal{M}$ can monotonically isotoped to another torus $\hat{M}$ that is $\epsilon$-close (in $C^0$ norm) to $\gamma$ \end{thm} 
Below we will call such $\hat{M}$ $\epsilon$-thick knots. Hence in many ways understanding (2-convex) tori is reduced to studying knots. As a corollary, we then show (analogous to theorem 1.1 above): 
\begin{thm} $\textbf{(Hypersurface to skeleton)}$ Given any 2-convex closed hypersurface $\mathcal{M}^n \subset \R^{n+1}$, there exists a skeleton $\gamma$ so that for all $\epsilon > 0$ $\mathcal{M}$ can be monotonically isotoped to an $\epsilon$ thick skeleton $\gamma$ in $\R^{n+1}$ (in the same sense as $\epsilon$ thick knot above). 
\end{thm} 
By skeleton here we mean a (very possibly nonsmooth) set of two types:
\begin{enumerate}
\item A point
\item consisting of embedded $S^1$'s connected by (individually) embedded, possibly intersecting intervals like below:
\end{enumerate}
 \begin{center}
$\includegraphics[scale = .27]{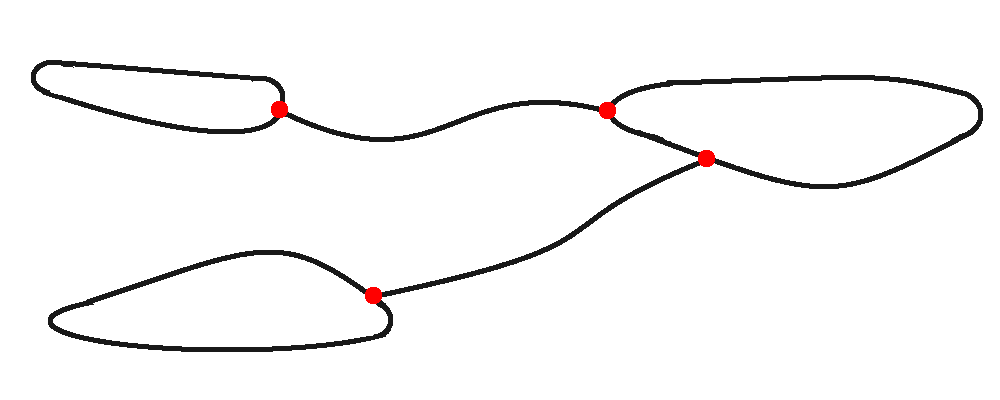}$
\end{center}
It might be interesting to point out that the these two statements rely on Theorem 1.1 above (that is, just from the surgery procedure). A further corollary of this is an aforementioned finiteness theorem for 2-convex hypersurfaces as an `extrinsic" sort of Cheeger's compactness theorem \cite{Che}, up to isotopy:
\begin{cor} $\textbf{(Finiteness Theorem})$ Let $\Sigma(d, C, \alpha)$ be the set of closed 2-convex hypersurfaces $\mathcal{M}$ with diam$(\mathcal{M}) < d$ (or equivalently up to translation, $\mathcal{M} \subset B_d(0)$), $H < C$ (initially), and is (2-sided) $\alpha$ non-collapsed for a fixed $\alpha > 0$. Then $\Sigma(d, C, \alpha)$ up to isotopy consists of finitely many hypersurfaces, in fact at most $2^{2n} \frac{( 12dC\sqrt{n})^n}{\alpha^n}$. 
\end{cor}
As mentioned above, at the time this manuscript was being written results similar to theorem 1.3 and 1.4 were attained using similar surgery techniques by Buzano, Halshofer, and Hershkovits in \cite{BHH} and \cite{BHH1}. Their isotopies aren't monotonic which complicates some of the arguments quite a bit but they notably manage to preserve 2-convexity. For complete fidelity we state their theorems as stated in their papers: 
\begin{thm}(Main Theorem in \cite{BHH}) The moduli space of 2-convex embedded $n$-spheres in $\R^{n+1}$ is path connected for every dimension $n$, i.e. $\pi_0(\mathcal{M}_n^{2-conv}) = 0$. \end{thm}

\begin{thm}(Theorem 1.3 in \cite{BHH1}) The path components of the moduli space of two-convex embedded tori are given by 
\begin{center}
$\pi_0(\mathcal{M}^{2-conv}(S^{n-1} \times S^1)) \cong \mathcal{K}(\mathcal{M}^{2-conv}(S^1 \times S^1))$ when $n=2$, and $\cong 0$ when $n \geq 3$.
\end{center}
where $\mathcal{K}$ denotes the set of knot classes. This means that $\mathcal{M}^{2-conv}(S^{n-1} \times S^1)$ is path connected for $n \geq 3$, while for $n=2$ we have that two mean-convex embedded tori are in the same path-component of their moduli space if and only if they have the same knot class.
\end{thm}

First though because of how central to our argument the flow with surgery is to our arguments we provide a brief account of it and related definitions; afterwards we explain the results:
\tableofcontents

\subsection{Acknowledgements} 
The author wishes to thank his advisor Rick Schoen, in addition to his great patience and generosity, for mentioning the Smale conjecture and suggesting flow methods might eventually be used to reprove it, which inspired this work (although the end product falls far short of the goal). The author also thanks Hung Tran and Dave Wiygul for several stimulating conversations throughout the course of this work.

\section{A brief introduction to the mean curvature flow with surgery}
\subsection{Classical formulation of the mean curvature flow}
There are multiple ways to define the mean curvature flow of a subset of $\R^{n+1}$ (Brakke flow, level set flow, etc.) but in this article we will focus on the differential geometric, or ``classical," definition of mean curvature flow for smooth embedded hypersurfaces of $\R^{n+1}$ (for a nice introduction, see \cite{Mant}). Let $\mathcal{M}$ be an $n$ dimensional manifold and let $F: \mathcal{M} \to \R^{n+1}$ be an embedding of $\mathcal{M}$ realizing it as a smooth closed hypersurface of Euclidean space - which by abuse of notation we also refer to $\mathcal{M}$. Then the mean curvature flow of $\mathcal{M}$ is given by $\hat{F}: \mathcal{M} \times [0,T) \to \R^{n+1}$ satisfying (where $\nu$ is outward pointing normal and $H$ is the mean curvature):
\begin{equation}
\frac{d\hat{F}}{dt} = -H \nu, \text{ } \hat{F}(\mathcal{M}, 0) = F(\mathcal{M})
\end{equation} 
(It follows from the Jordan seperation theorem that closed embedded hypersurfaces are oriented). Denote $\hat{F}(\cdot, t) = \hat{F}_t$, and further denote by $\mathcal{M}_t$ the image of $\hat{F}_t$ (so $\mathcal{M}_0 = \mathcal{M}$). It turns out that (2.1) is a degenerate parabolic system of equations so take some work to show short term existence (to see its degenerate, any tangential perturbation of $F$ is a mean curvature flow). More specifically, where $g$ is the induced metric on $\mathcal{M}$:
\begin{equation}
 \Delta_g F = g^{ij}(\frac{\partial^2 F}{\partial x^i \partial x^j} - \Gamma_{ij}^k \frac{\partial F}{\partial x^k}) = g^{ij} h_{ij} \nu = H\nu
\end{equation}
Now one could apply for example deTurck's trick to reduce the problem to a nondegenerate parabolic PDE (see for example chapter 3 of \cite{Bake}) or similarly reduce the problem to an easier PDE by writing $\mathcal{M}$ as a graph over a reference manifold (see \cite{Mant}). At any rate, we have short term existence. Let us record associated evolution equations for some of the usual geometric quantities: 
\begin{itemize}
\item $\frac{\partial}{\partial t} g_{ij} = - 2H h_{ij}$
$\medskip$

\item $\frac{\partial}{\partial t} d\mu = -H^2 d\mu$
$\medskip$

\item $\frac{\partial}{\partial t} h^i_j = \Delta h^i_j + |A|^2 h^i_j$
$\medskip$

\item $\frac{\partial}{\partial t} H = \Delta H + |A|^2 H$
$\medskip$

\item $\frac{\partial}{\partial t} |A|^2 = \Delta |A|^2 - 2|\nabla A|^2  + 2|A|^4$
\end{itemize}
So, for example, from the heat equation for $H$ one sees by the maximum principle that if $H > 0$ initially it remains so under the flow. There is also a more complicated tensor maximum principle by Hamilton originally developed for the Ricci flow (see \cite{Ham1}) that says essentially that if $\mathcal{M}$ is a compact manifold one has the following evolution equation for a tensor $M$: 
\begin{equation}
\frac{\partial M}{\partial t} = \Delta M + \Phi(M)
\end{equation} 
and if $M$ belongs to a convex cone, then if solutions to the system of ODE
\begin{equation}
\frac{\partial M}{\partial t} = \Phi(M)
\end{equation}
stay in that cone then solutions to the PDE (2.2) stay in the cone too (essentially this is because $\Delta$ ``averages"). So, for example, one can see then that convex surfaces stay convex under the flow very easily this way using the evolution equation above for the Weingarten operator. Similarly one can see that $\textbf{2-convex hypersurface}$ (i.e. for the two smallest principal curvatures $\lambda_1, \lambda_2$, $\lambda_1 + \lambda_2 > 0$ everywhere) remain 2-convex under the flow.
$\medskip$

Another important curvature condition in this paper is $\alpha \textbf{non-collapsing}$: a mean convex hypersurface $M$ is said to be 2-sided $\alpha$ non-collapsed for some $\alpha > 0$ if at every point $p\in M$, there is an interior and exterior ball of radius $\alpha/H(p)$ touching $M$ precisely at $p$. This condition is used in the formulation of the finiteness theorem. It was shown by Ben Andrews in \cite{BA} to be preserved under the flow (a sharp version of this statement, first shown by Brendle in \cite{Binsc} and later Haslhofer and Kleiner in \cite{HK3}, is important in \cite{BH} when MCF+surgery to $n=2$ was first accomplished). 
$\medskip$

Finally, perhaps the most important manifestation of the maximum principle is that if two compact hypersurfaces are disjoint initially they remain so under the flow. So, by putting a large hypersphere around $\mathcal{M}$ and noting under the mean curvature flow that such a sphere collapses to a point in finite time, the flow of $\mathcal{M}$ must not be defined past a certain time either in that as $t \to T$, $\mathcal{M}_t$ converge to a set that isn't a manifold.  Note this implies as $t \to T$ that $|A|^2 \to \infty$ at a sequence of points on $\mathcal{M}_t$; if not then we could use curvature bounds to attain a smooth limit $\mathcal{M}_T$ which we can then flow further, contradicting our choice of $T$. 
$\medskip$

So, to use mean curvature flow to study compact hypersurfaces one is faced with finding a way to extend the flow through singularities. To do this, one could define the mean curvature flow on more general sets than manifolds (cf. the Brakke flow, level set flow) but these have the disadvantage of moving outside the realm of smooth differential geometry, contrary to our goal.
$\medskip$

Now as we mentioned before, for a 2-convex hypersurface $\mathcal{M}$ regions on $\mathcal{M}_t$ that develop high curvature can be basically classified due to original curvature assumptions (namely 2-convexity) whose topology we are able to easily identify, and so by defining a surgery theory to cut these pieces out and glue in caps smoothly one can stay within differential geometry and also have good understanding of the topology of the flow past the first singular time. However of course there are a great many technical details to overcome to make this precise and there are a couple somewhat different schemes put forth; one essentially by Huisken, Sinestrari, and Brendle (see \cite{HS2}, \cite{BH}) and another newer way by Haslhofer and Kleiner (see \cite{HK1}). 
$\medskip$

Here as mentioned above we follow the scheme set out by Huisken and Sinestrari in \cite{HS2} (which develops mean curvature flow with surgery for hypersurfaces in $\R^{n+1}$, where $n \geq 3$), the scheme which was later extended by Brendle and Huisken's later additions in \cite{BH} for the $n=2$ case. For the sake of brevity we only state without proof what we need in the sequel - we are necessarily going to elide some details as a result and the reader is encouraged to peruse the original sources.
\subsection{Mean curvature flow with surgery according to Huisken, Sinestrari, and Brendle}
$\medskip$ 

Amongst many constants in rigourously defining the mean curvature flow with surgery, there are three constants most relevant to our needs which we denote $H_1 < H_2 <  H_3$. The mean curvature flow with surgery is defined by taking a compact (2-convex) hypersurface $\mathcal{M}$ and to start flowing it by the mean curvature flow. It turns out with appropriate choice of constants (including of course the $H_i)$ that when there is a point $p$ where $H = H_3$ on $\mathcal{M}_t$ there are two possibilities. The first is that $\mathcal{M}_t$ will be uniformly convex (in that all the principle curvatures $\lambda_i$ are bounded below by some positive constant). Note in this case $\mathcal{M}_t$ is hence star shaped so $M_t$ is diffeomorphic to $S^n$. The other possibility is that in a neighborhood of $p$ $\mathcal{M}_t$ will locally look like $S^{n-1} \times I$ (that is after rescaling, be a graph of small norm (say $\epsilon > 0$ over a cylinder in $C^k$ norm ($k >2$, $k$ depending on fixed parameters). if near $p$ the manifold is locally cylindrical (or locally like a neck, if you will), there is a $\textbf{neck continuation theorem}$ that roughly says that $p$ is in fact a point in a maximally extended neck $\mathcal{N}$, a \textbf{hypersurface neck} region being a region $\mathcal{N} \subset \mathcal{M}$ is where every near every point $p \in \mathcal{N}$, $\mathcal{N}$ is locally graphical over a cylinder \footnote{there are actually several definitions of neck one wants to use in completely defining the surgery that one shows are equivalent; for examples a curvature neck is a region where the second fundamental form is close to that of a cylinder (after rescaling) - for more details see section 3 of \cite{HS2} - for our purposes this simple definition suffices}. There are several cases to consider for $\mathcal{N}$: 
\begin{center}
\begin{enumerate}
\item it is bordered on both ends by regions of low curvature $H \sim H_1$
\item it is bordered on one or both ends by caps (one of the ends could be bordered by a region of low curvature)
\item both ends of the neck "meet" to form a loop
\end{enumerate}
\end{center}
If regions of low curvature border the neck, the neck is cut and caps are inserted where $H \sim H_2$, separating regions where $H \leq H_2$ and high curvature regions that belong to the neck. These regions are diffeomorphic to either $S^n$ or $S^{n-1} \times S^1$ as the figure below illustrates
  \begin{center}
$\includegraphics[scale = .3]{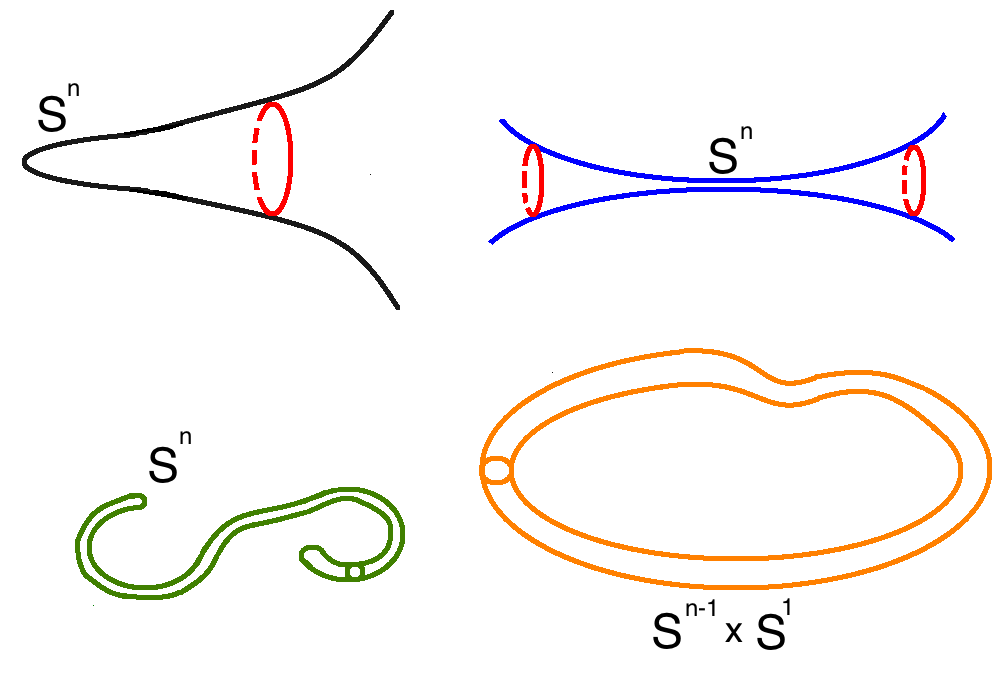}$
\end{center}
Its important to note that the caps inserted from the surgery process are convex (caps/ends of a neck can also be from $\mathcal{M}_t$ ``naturally" tapering off and in fact will also be convex. We will say (using a non-standard term) that a $\textbf{standard surgery}$\footnote{to differentiate from the degenerate cases where the whole manifold is high curvature} is performed at $p$ if $\mathcal{M}$ is cut at $p$ and caps are glued in (disconnecting $\mathcal{M}$ locally). So if the hypersurface neck is bounded by two low curvature regions, there will be two standard surgeries along the neck 
 \begin{center}
$\includegraphics[scale = .3]{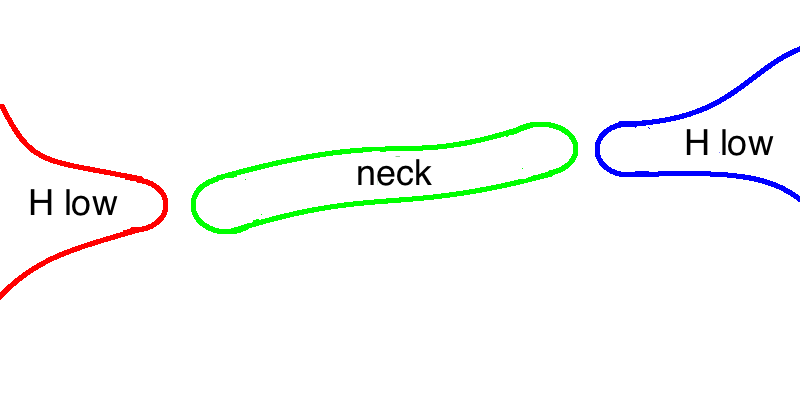}$

\end{center}
Continuing on, the high curvature components in the surgery algorithm are thus classified topologically and thus deleted, and the flow is continued on the left over low curvature regions (if any) until a point where $H = H_3$ is encountered again and the process is repeated.  
$\medskip$

Finally, since $H \leq H_3$ on the high curvature regions, at every surgery time there is a positive lower bound on the volume of $\mathcal{M}_t$ removed by surgery. Since mean curvature flow decreases volume (in fact, $\frac{d}{dt} vol(\mathcal{M}_t) = -\int_{\mathcal{M}_t} H^2 d\mu_t$) there are a finite number of surgeries until $\mathcal{M}_t$ is exhausted (i.e. there are no low curvature regions left over). This essentially shows Theorem 1.1 above. Also, important for us is that fixing all choices of parameters, the number of surgeries is well defined.

\section{Path to round sphere}

In this section we wish to prove theorem 1.3, namely that any 2-convex simply connected hypersurface $\mathcal{M}^n \subset \R^{n+1}$ is isotopic to the standard round sphere through a monotone isotopy. To do this of course we will use the mean curvature flow with surgery. Towards this end, we $\textbf{induct on the number of standard surgeries}$ encountered along the flow with surgeries $\mathcal{M}_t$ as described in the previous sections \footnote{note that, upon fixing parameters, the number of surgeries is well defined}. Let $T_0$ be the first surgery time of $\mathcal{M}_t$ (of course, by the assumption of mean convexity, the flow gives a monotone isotopy from $\mathcal{M}$ to $\mathcal{M}_{T_0}$). 
$\medskip$

The base case in the induction argument is if there are no standard surgeries encountered along $\mathcal{M}_t$. At time $T_0$ then, somewhere on $\mathcal{M}_{T_0}$ we have $H = H_3$. $\mathcal{M}_{T_0}$ either contains a hypersurface neck or it doesn't. If it doesn't, then as mentioned in the previous section $\mathcal{M}_{T_0}$ is uniformly convex. If there are necks, then from neck continuation theorem they can be continued until they hit a low curvature region or are capped off (furthermore, these maximally extended necks are disjoint). If there were low curvature regions, we would cut the neck out and glue in surgery caps (leaving the low curvature regions behind) contrary to our assumption of no standard surgeries. Hence since $\mathcal{M}_{T_0}$ is connected it is either  ``no neck" or ``all neck;" see the picture below for the two cases:
 \begin{center}
$\includegraphics[scale = .3]{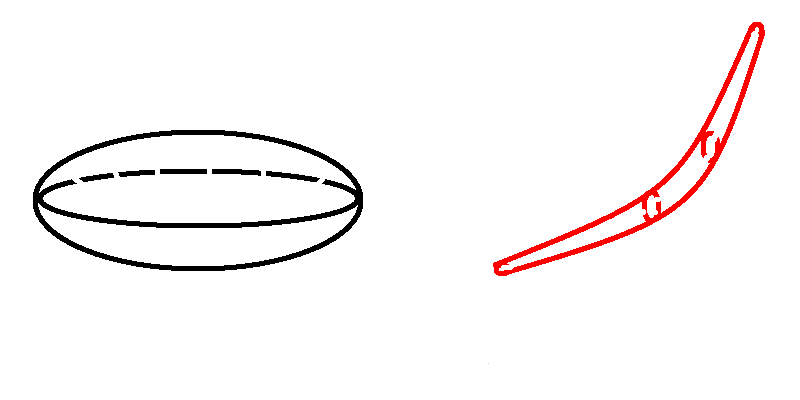}$
\smallskip

$\textbf{Figure 1}$
\end{center}
$\medskip$

If $\mathcal{M}_{T_0}$ is strictly convex as in the first case, then namely it is starshaped with respect to any point in the interior of the domain it bounds. Hence by picking such a point and taking a small sphere $S$ around it, we have a straight line homotopy of $\mathcal{M}$ to $S$ which of course is also monotone. 
 \begin{center}
$\includegraphics[scale = .3]{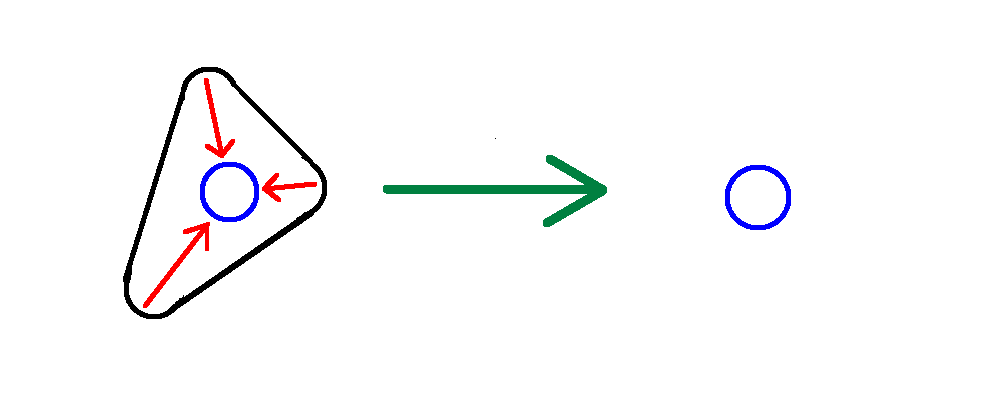}$
\smallskip

$\textbf{Figure 2}$
\end{center}
$\medskip$

 If $\mathcal{M}_{T_0}$ is as in the second case though things are more complicated and require understanding the caps better. Now, in our particular case ($T_0$ the first surgery time) we know that there wern't any prior surgeries, but because we rely on induction we must allow for the presence of prior surgeries along the flow. In particular we will show:
 \begin{prop}If $\mathcal{M}_{T_0}$ coincides with a maximal hypersurface neck that ends in two caps in the sense of the neck continuation theorem, then $\mathcal{M}_{T_0}$ is isotopic to a round sphere through a monotone isotopy. 
 \end{prop}
 \begin{pf} The basic idea is to first isotope (monotonically) $\mathcal{M}_{T_0}$ to a more explicitly described hypersurface neck with caps whoose isotopy to $S^n$ is easy to describe. There are two types of points along our surface in this case, a point at the center of an neck (for some $\epsilon,k$ - these can be controlled by changing the surgery parameters) or a point on a cap; the caps are given by gluing a disc along a cross section the boundary of the maximal neck. Let us denote $\mathcal{M}_{T_0} := \mathcal{N} \cup (D_1 \cup D_2)$ to reflect this. 
$\medskip$

We recall every point of a hypersurface neck $\mathcal{N}$  can be written locally as the graph over a cylinder. Since $\mathcal{M}_{T_0}$ is compact we can pick a finite set of points $\{p_i\}$ with neighborhoods $U_i$ such that $\{p_i, U_i\}$ cover $\mathcal{N}$. Let $\gamma_i$ be the principal axis of these cylinders and let $\gamma$ be these $\gamma_i$ glued together. Let us further denote $\hat{\gamma}$ to be a small uniform tubular neighborhood of $\gamma$, taken small enough so that $\hat{\gamma}$ lies in the interior of $\mathcal{N}$.
$\medskip$

Then $\mathcal{N}$ is a graph over $\hat{\gamma}$, and we could use a straight line homotopy to provide a monotone isotopy from one to the other except we also have to worry about the caps. 
$\medskip$

To deal with this we have to understand the caps. From the proof of the neck continuation theorem we know the caps are of one of two types; either the cap comes from a ``recent surgery" or not. If not, then it is shown that the cap is convex; suppose first our caps $D_1$ and $D_2$ are convex. Then we see that it is a graph over a hemisphere of $S^{n-1}$ situated along the boundary of $\mathcal{N}$ so that extending the straightline homotopy along the neck we get a monotone isotopy to $\hat{\gamma}$.
 \begin{center}
$\includegraphics[scale = .3]{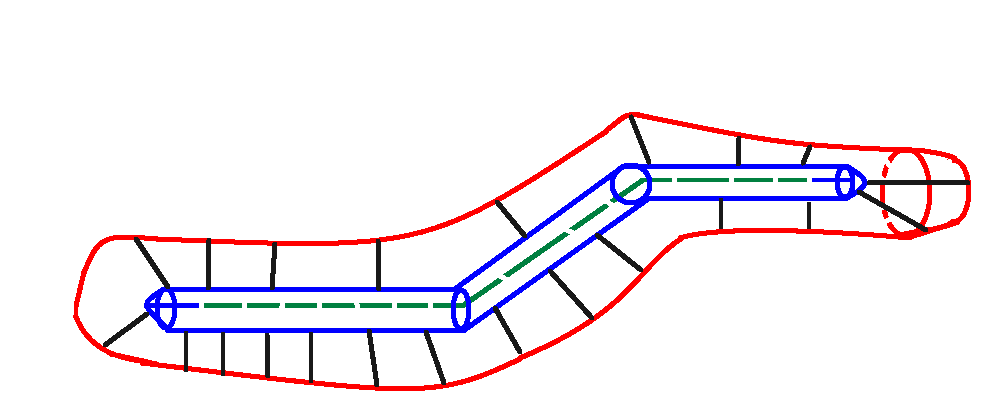}$
\smallskip

$\textbf{Figure 3}$
\end{center}
$\medskip$

Suppose though one of the caps comes from a recent surgery. Since the surgery caps are initially convex by design, at the time of the surgery (call it $T_c$) we then have that the cap is a graph over an appropriate hemisphere. We recall from \cite{HS2} that after a surgery is performed another surgery won't be performed nearby if it was recent (essentially, the curvature wouldn't have time to increase from $~H_2$ to $H_3$)\footnote{ see theorem 8.2 in \cite{HS2} - unfortunately it seems too much to fully explain this in the introduction to surgery}, so that flow with surgery is an isotopy from $T_c$ to $T_0$ in that region. Hence if we straight line isotope the cap at $T_c$ as below to a cross section of $\mathcal{N}$ by straightline homotopy, we can the deform $\mathcal{M}_{t}$ for all $t \in [T_c, T_0]$ so that at $T_0$ $\mathcal{M}_t$ (deformed) instead of having a cap coming from surgery of indeterminate geometry has as a cap a cross-section of the neck (of course, we can ``round" the edges on the cross section to preserve smoothness). We may isotope this then to $\hat{\gamma}$, monotonically, using the straightline homotopy from before. 
 \begin{center}
$\includegraphics[scale = .3]{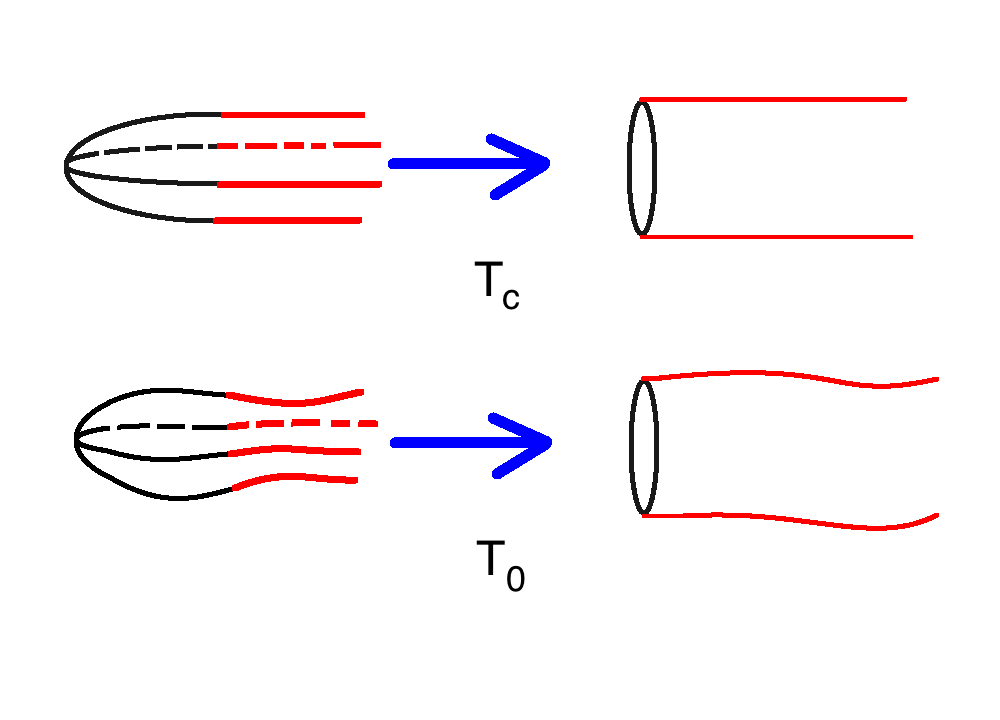}$
\smallskip

$\textbf{Figure 4}$
\end{center}
$\medskip$

So in any case we have a monotone isotopy of $\mathcal{M}_{T_0}$ to $\hat{\gamma}$. To monotonically isotope this down to the sphere we essentially push one end in; we proceed by induction on the number $N$ of straight segments $\gamma_i$ that $\gamma$ is composed of. If $N =1$, then we can write $\hat{\gamma}$ as two convex caps glued together on opposite sides to a cylinder, so the homotopy is just to retract the cylinder to $S^{n-1}$ by a straightline homotopy (since the caps are convex, it is starshaped). 
  \begin{center}
$\includegraphics[scale = .3]{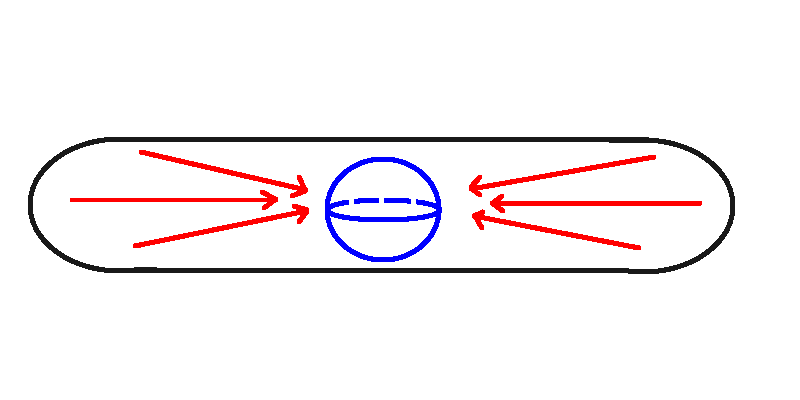}$
\smallskip

$\textbf{Figure 5}$
\end{center}

Now suppose the statement is known for $N = m -1$ and that $N = m$. Then we retract along the ``cylinder" of $\gamma_1$ by a straightline homotopy to the boundary of $\hat{\gamma_1}$ and $\hat{\gamma_2}$, and then we further homotope (of course, by straightline homotopy) to a convex cap of a neck composed of  $m-1$ straight segments and so conclude by induction
\begin{center}
$\includegraphics[scale = .3]{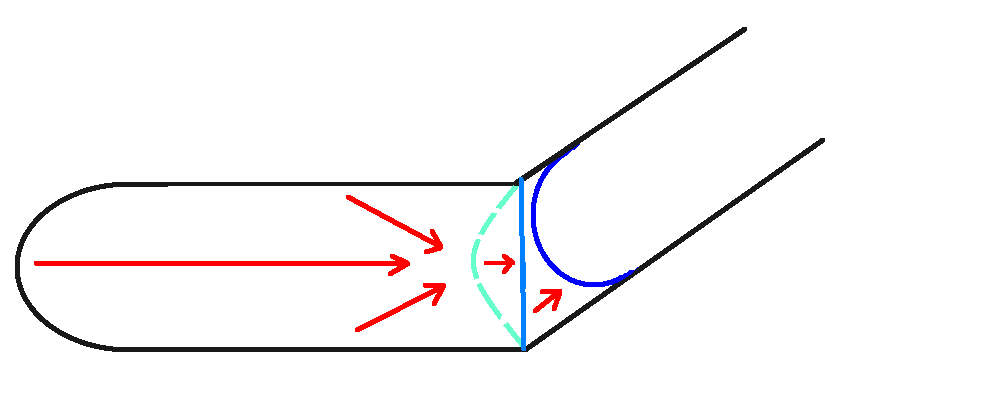}$
\smallskip

$\textbf{Figure 6}$
\end{center}
$\medskip$

 \end{pf}

Now suppose that $\mathcal{M}_{T_0}$ has at least one point where a standard surgery should take place by the surgery algorithm and consider the diagram below of $\mathcal{M}$ as a graph (with the $\mathcal{M}_i$ as vertices\footnote{if the $\mathcal{M}_i$ are adjacent along a neck that contains two standard surgery positions, one of the $M_i$ contains one of the surgery spots closest to it (but it doesn't matter which)}, and the red lines indicating approximately where a standard surgery should be done).
 \begin{center}
$\includegraphics[scale = .4]{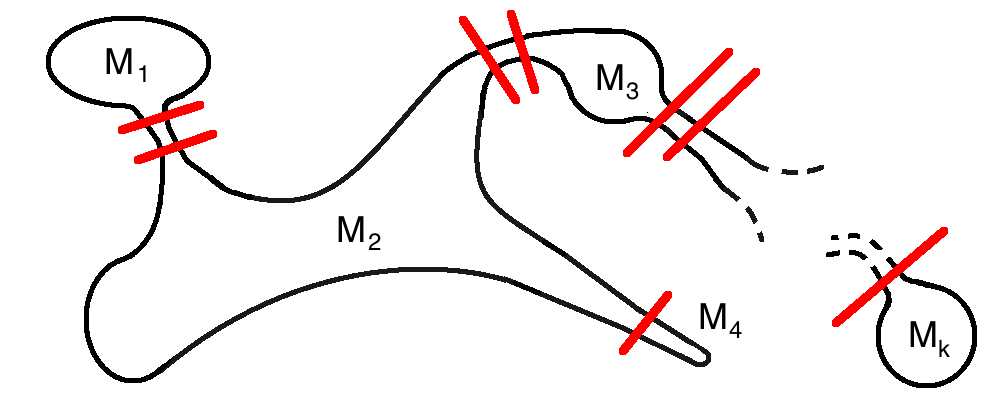}$
\smallskip

$\textbf{Figure 7}$ 
\end{center}
From the classification theorem 1.1 above, we see since $\mathcal{M}$ is simply connected that each of the $\mathcal{M}_i$ must be too (the only choices are connect sums of tori or spheres) and as a graph $\mathcal{M}$ is a tree (no cycles). Hence, picking a single neck to do surgery along will split $\mathcal{M}_{T_0}$ into two hypersurfaces $\mathcal{A}$ and $\mathcal{B}$ diffeomorphic to $S^n$ whose flow with surgeries (using the same parameters) have strictly less surgeries than $\mathcal{M}$ does. Hence it suffices to show:
\begin{prop} Suppose $\mathcal{M} = \mathcal{A} \# \mathcal{B}$ where $\mathcal{A}$ and $\mathcal{B}$ are two hypersurfaces diffeomorphic to $S^n$ joined together by a hypersurface neck $\mathcal{N}$ (in application, a maximal hypersurface neck found by the neck continuation theorem) such that 
\begin{enumerate}
\item Doing surgery at one point of $\mathcal{N}$ where $H \sim H_1$, as described in the surgery procedure to attain $\mathcal{M}^+$, leaves us with $\mathcal{A}$ and $\mathcal{B}$. 
\item $\mathcal{A}$ and $\mathcal{B}$ are both monotonically isotopic to round spheres.
\item There exists open sets $U_1, U_2$ such that $\mathcal{A} \subset U_1$ and $\mathcal{B} \subset U_2$, where
\begin{enumerate}
\item $U_1$ and $U_2$ are disjoint 
\item there is an $\epsilon > 0$ such that $\mathcal{M}^+ \cap U_1$, $\mathcal{M}^+ \cap U_2$ have $\epsilon$-tubular neighborhoods contained in $U_1$, $U_2$ respectively.
\end{enumerate}
\end{enumerate}
Then $\mathcal{M}$ is monotonically isotopic to a round sphere.
\end{prop}
\begin{pf}
First note that, because after cutting and pasting the result caps are a fixed distance apart and that each of the $\mathcal{M}_I$ above are disjoint compact hypersurfaces, we are indeed in a situation that satisfies (3) above. The situation is displayed abstractly in the diagram below ($\mathcal{A}$ and $\mathcal{B}$ might be very complicated): 
 \begin{center}
$\includegraphics[scale = .3]{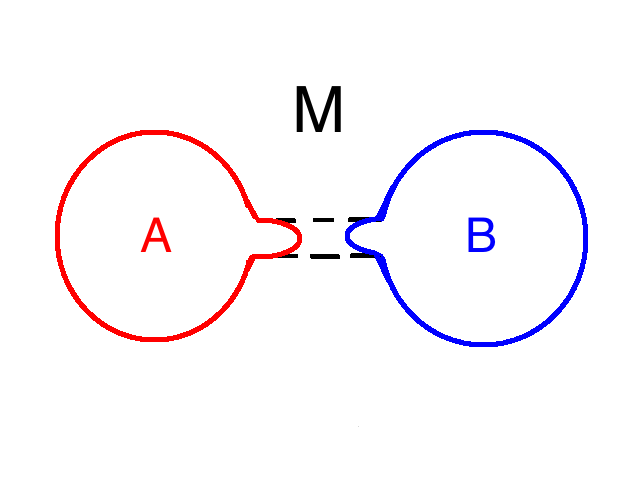}$
\smallskip

$\textbf{Figure 8}$ 
\end{center}
Denote by $H_\mathcal{A}: S^{n} \times [0,1] \to \R^{n+1}$ and $H_\mathcal{B}: S^{n} \times [0,1] \to \R^{n+1}$ the homotopies bringing $\mathcal{A}$ and $\mathcal{B}$ to round spheres; furthermore denote their time slices $\mathcal{A}_t$ and $\mathcal{B}_t$ respectively. Then the plan is to ``patch in" the neck between $\mathcal{A}_t$ and $\mathcal{B}_t$ in a continuous (in time), embedded fashion so that we have an isotopy $\mathcal{A}_t \#_t \mathcal{B}_t$ from $\mathcal{M}$ to round spheres $\mathcal{A}_1$, $\mathcal{B}_1$ connected by a thin neck. Then we will homotope this to a round sphere as in the picture below:
 \begin{center}
$\includegraphics[scale = .3]{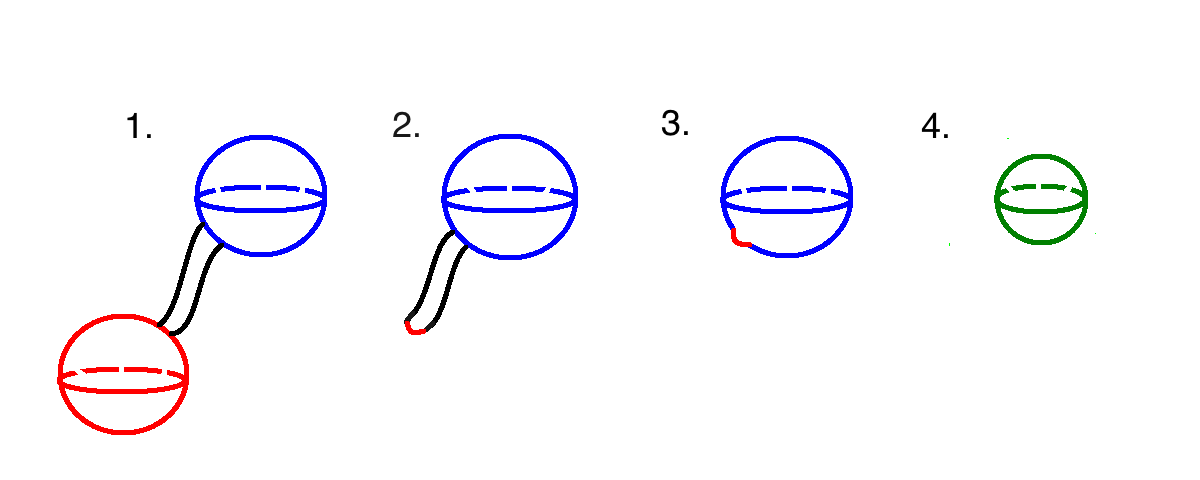}$
\smallskip

$\textbf{Figure 9}$ 
\end{center}
Some care has to be taken because we want all these homotopies to be embeddings at each time slice though; first off it is conceivable that at some time $t_0$ (before we even get to the setup of the picture above), $\mathcal{A}_{t_0}$ and $\mathcal{B}_{t_0}$ intersect. For us though this is taken care of by monotonicity of the isotopies; in fact it implies $\mathcal{A}_t \subset U_1$, $\mathcal{B}_t \subset U_2$ for all $t$. We also have to define precisely how to ``extend the neck" in a way that will give a monotonic isotopy of $\mathcal{M}$ - again monotonicity of $\mathcal{A}$ and $\mathcal{B}$ lets us relatively easily construct a neck (hypersurface that can locally be written as a graph over a cylinder) that does the trick as we'll see below. 
$\medskip$

Now to begin, from smoothness of isotopy and compactness of $[0,1]$ there is a uniform upper bound on $|A|^2$ for $\mathcal{A}_t$ and $\mathcal{B}_t$ so that there is an $\eta >0$ where in any $\eta$ ball $B_P(\eta)$ of any point $P$ on $\mathcal{A}_t$ or $\mathcal{B}_t$, the hypersurfaces can be locally written as a graph over $T_P \{\mathcal{A}_t, \mathcal{B}_t, \}$. Furthermore, for the sequel we take $\eta$ small enough so that for any $P_1$ in $B_P(\eta)$, $|\nu(P) - \nu(P_1)| < \epsilon/1000$ say ($\epsilon$ as in the assumption).
$\medskip$

With that in mind, we homotope $\mathcal{M}$ near each of the places where caps would be inserted (this is preconditioning if you will - we aren't isotoping $\mathcal{A}$ or $\mathcal{B}$ yet). The picture near $\mathcal{A}$:
\begin{center}
$\includegraphics[scale = .28]{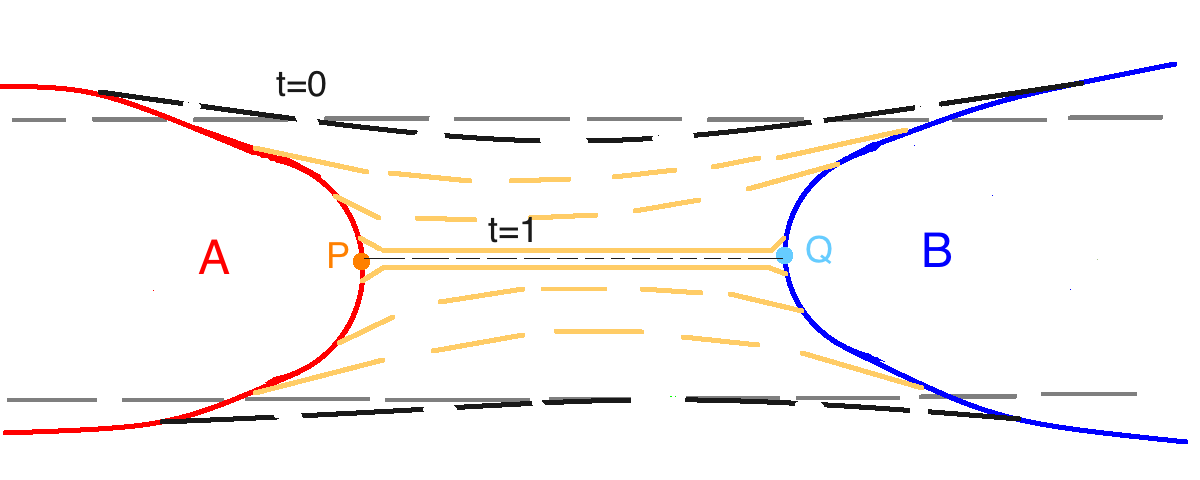}$

$\textbf{Figure 10}$ 
\end{center}
We are just pinching the neck and such a homotopy is clearly possible by a monotone homotopy since the convex caps put in place after a standard surgery are disjoint. More precisely choose this homotopy so that at $t=1$ above the deformed neck is precisely a round cylinder with diameter $d \leq \eta/4$ glued to $\mathcal{A}_0$, $\mathcal{N}_1$, which at the interface is most likely nonsmooth:
\begin{center}
$\includegraphics[scale = .35]{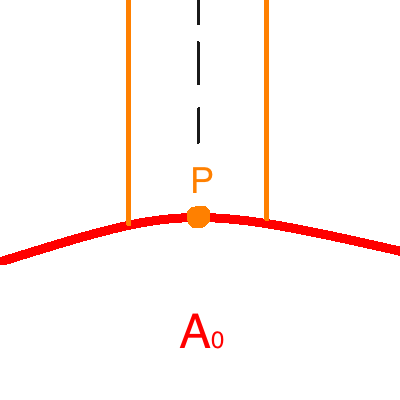}$

$\textbf{Figure 11}$ 
\end{center}

Here the central axis of this cylinder passes through (on the diagram) points $P \in \mathcal{A}$ and $Q \in \mathcal{B}$ which are at the centers of the caps so that the central axis is normal to $\mathcal{A}$ and $\mathcal{B}$ so that we are in the case below \begin{center}
$\includegraphics[scale = .3]{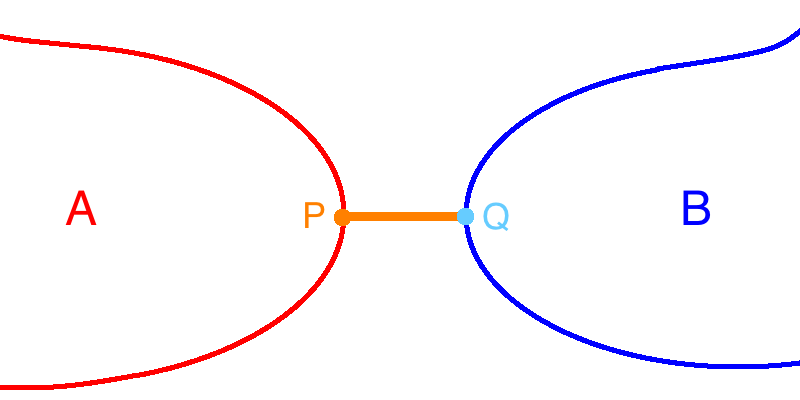}$
\smallskip

$\textbf{Figure 12}$ 
\end{center}
Now we start to let $\mathcal{A}$ and $\mathcal{B}$ isotope to round spheres - to define the connecting neck for positive times we keep track of $P(t)$ and $Q(t)$; here by $P(t)$ and $Q(t)$ we mean specifically with respect to the ``normalized" isotopy with derivative only in the normal component of the hypersurface (tangential perturbations don't actually affect the realization of the hypersurface in $\R^{n+1}$)\footnote{this is in a sense the same problem as the degeneracy of the mean curvature flow discussed in the introduction to the mean curvature flow equation}. We can focus our attention on extending the neck along the isotopy of $\mathcal{A}_t$ and follow a similar procedure for extending it along the isotopy of $\mathcal{B}_t$ 
 $\medskip$
 
Assuming without loss of generality then that our isotopy is ``normalized" as discussed in the previous paragraph and because the isotopy $\mathcal{A}_t$ is monotone and smooth in time, the path $\gamma_\mathcal{A}(t)$ of $P(t)$ in $\R^{n+1}$ is 
 \begin{enumerate}
 \item smooth in time.
 \item embedded,
 \item  $\nu(P(t)) \perp T_{P(t)} \mathcal{A}_t$, and
 \item  disjoint from $A_s$ for $s<t$
 \end{enumerate}
 \begin{center}
$\includegraphics[scale = .3]{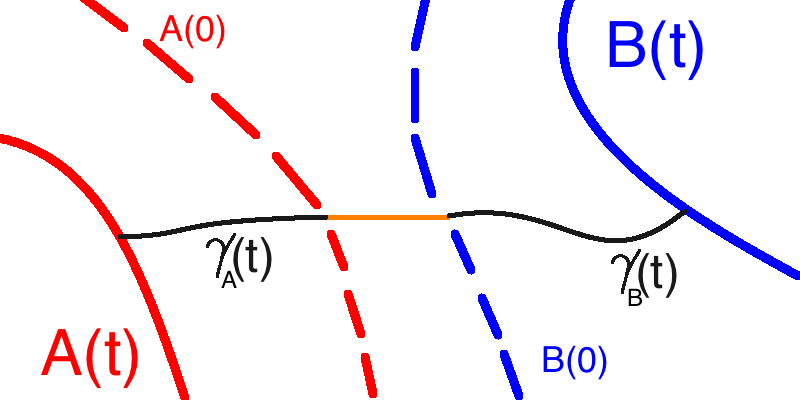}$
\smallskip

$\textbf{Figure 13}$ 
\end{center}
Of course, for a fixed time this isn't a neck but we see that its $d$-tubular neighborhood $\gamma_d$ ($d$ from above) is and we see (repeating the construction on $\mathcal{B}_t$) for each time $t$ we get apriori a nonsmooth (yet continuous) monotone isotopy $\mathcal{A}_t \#_t \mathcal{B}_t$; to deal with this we will redefine the isotopy by smoothing at the interface. 
$\medskip$

Now note for an isotopy $\mathcal{M}_t$ for say $t \in [\alpha, \beta] \subset [0,1]$,  $\mathcal{M}_t$ being a monotone isotopy is equivalent to, writing $\mathcal{M}_t$ locally as a graph of a function $f(x,t)$ over some reference manifold $\tilde{M}$ contained inside the hull of $\mathcal{M}_t$ for $t \in [\alpha, \beta]$ (using the convention that the normal vector is taken to be outward pointing), that for $s > t$, $f(x,s) < f(x,t)$. Since molification is linear and the molification of a positive function remains positive (mollifying use a positive bump function) we see that monotonicity is thus preserved. 
$\medskip$

We see for our particular case that we can't expect to find a reference manifold that all of $\mathcal{A}_t \#_t \mathcal{B}_t$ can be written as a graph over that works for all $t \in [0,1]$; this is bad news because although certainly we can split $[0,1]$ into subintervals where we can do this, it could concievably be non-monotone or even noncontinuous when we transition from one reference manifold to the next, at least if not done carefully. Luckily, however, the only regions of $\mathcal{A}_t \#_t \mathcal{B}_t$ that need to be smoothed out are the interfaces, so if we can find a smooth hypersurface $\mathcal{K}$ such that 
\begin{enumerate}
\item  $\mathcal{K}$ is contained in the hull of $\mathcal{A}_t \#_t \mathcal{B}_t$ for all $t \in [0,1]$,
\item For all $t \in [0,1]$ the neck interface of $\mathcal{A}_t \#_t \mathcal{B}_t$ can be written locally as a graph over $\mathcal{K}$, and
\item away from the interface, the height function (function corresponding to graph) of the region of $\mathcal{A}_t \#_t \mathcal{B}_t$ that is graphical over $\mathcal{K}$ is smooth 
\end{enumerate}
We can then do our mollification trick to get a family of smooth hypersurfaces $\overline{\mathcal{A}_t \#_t \mathcal{B}_t}$ that is monotone in time (note by assumption (3b), they are embedded). Furthermore, since the isotopies $\mathcal{A}_t$ and $\mathcal{B}_t$ are smooth, the interface varies continuously in time; hence we can use a smoothly varying family of molifiers so that $\overline{\mathcal{A}_t \#_t \mathcal{B}_t}$ varies smoothly in time, so that the family is indeed an isotopy. 
$\medskip$

A natural candidate for such a $\mathcal{K}$ is $\gamma_{d/2}(1)$ (or some other small tubular neighborhood of $\gamma(1)$). We see from the construction of our ``rough" isotopy that it is contained in $\mathcal{A}_t \#_t \mathcal{B}_t$ for all $t \in [0,1]$; however it is possible that the interface is not a graph over it! Here is what the picture could be at the interface with $\mathcal{A}_t$ for example:
\begin{center}
$\includegraphics[scale = .4]{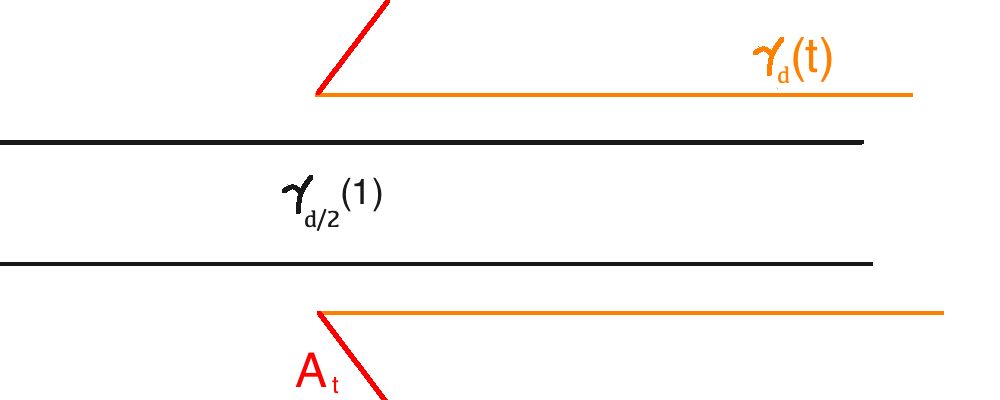}$

$\textbf{Figure 14}$ 
\end{center}
The problem being, more precisely, that (in a neighborhood of) the interface, the normal $\nu(x)$ of $\mathcal{A}_t$ (and/or $\mathcal{B}_t$) might be so that $\langle \nu(x), \nu(P(t)) \rangle \leq 0$ at some points $x \in A(t)$ 
$\medskip$

To overcome problems like this, we should perturb $\mathcal{A}_t$ and $\mathcal{B}_t$ near $P(t)$ and $Q(t)$ (name the perturbed points the same) so that for all $t$ in fact $\langle \nu(x), \nu(P(t)) \rangle > 0$ in a uniform neighborhood of the interface. The perturbation must be done in a smooth, monotonic way. 
$\medskip$

To see we can do this, we recall we choose $\eta$ so small that for $P_1 \in B_\eta(P)$, $P \in \mathcal{A}_t$ or $\mathcal{B}_t$, that $|\nu(P_1) - \nu(P)| < \epsilon/1000$ and further choose $d < \eta/4$, so if we add a smooth positive bump function (writing $\mathcal{A}_t$, $\mathcal{B}_t$ locally as a graph) at every time the to obtain hypersurfaces $\tilde{\mathcal{A}_t}$ and $\tilde{\mathcal{B}_t}$, and keep $\gamma(t), \gamma_d(t)$ as before, the interface will be a graph over $\mathcal{A}_t \#_t \mathcal{B}_t$ for each time slice (note we can take the perturbation small enough so that $\tilde{\mathcal{A}_t}$ and $\tilde{\mathcal{B}_t}$ will be embedded (i.e. not intersect) by assumption (3b)). 
\begin{center}
$\includegraphics[scale = .4]{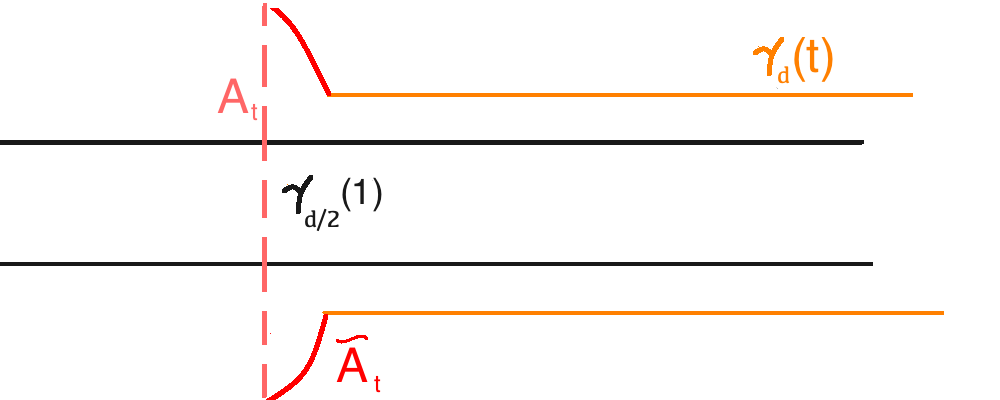}$

$\textbf{Figure 15}$ 
\end{center}
Since $\mathcal{A}_t$ and $\mathcal{B}_t$ are smooth and isotopic, we can arrange so that these perturbed families $\tilde{\mathcal{A}_t}$ and $\tilde{\mathcal{B}_t}$ are in fact monotone isotopies in their own right, in that they (as a family in $t$) vary smoothly in time and are monotone. In more detail, note since $\mathcal{A}_t$ and $\mathcal{B}_t$ vary smoothly in time, there exists a $\delta$ so that for, for any point $P \in \mathcal{A}_s$ or $\mathcal{P}_s$, the isotopy can be written locally as an evolution of graphs $f(x,t)$ for say $|t - s| < \delta$ over some hyperplane $H_P$ (say a translation of the tangent plane of $P$), where $f$ is smooth in both $x$ and $t$. 
$\medskip$

Since the isotopies are monotone, $f(x,t)$ vary monotonically in that as described before for $t_1 > t$, $f(x, t_1) < f(x,t)$, so if we add a fixed positive bump function $\chi$  we get for small time a perturbed hypersurface that varies smoothly in time and is monotone. 
$\medskip$

So, first we perturb $\mathcal{A}_t$ and $\mathcal{B}_t$ for $t \in [0, \delta/2)$, and then we move the hyperplane in advance to a hyperplane $\overline{H}_P$, so that locally at $P$ the hypersurface (either $\mathcal{A}$ or $\mathcal{B})$ can locally be written as a graph of $\overline{f}(x,t)$ over $\overline{H}_P$ for say $t \in [ \delta/2, 3\delta/2)$; at $t = \delta/2$ the hypersurface can either be written locally near $P(t)$ as a graph over $H_P$ or $\overline{H}_P$. 
$\medskip$

It remains then to find $\overline{\chi}$ so that $f(x,\delta/2) + \chi(x) = \overline{f}(y, \delta/2) + \overline{\chi}(y)$, where $x$ are the coordinates of $H_P$ and $y$ are the coordinates on $\overline{H}_P$. Note though that $y$ is related to $x$ by a rotation and translation; $x = R(y) + \mathbf{C}$. Using this we see we $\overline{\chi}(y) 
= f(Ry + C, \delta/2) - \overline{f}(y,\delta/2) + \chi(Ry + \mathbf{C})$. From here we see how to perturb $\mathcal{A}_t$ and $\mathcal{B}_t$ for all $t \in [0,1]$ with the desired properties. 
$\medskip$

Hence (again, using the same $\gamma, \gamma_d$ as before) $\tilde{\mathcal{A}_t} \#_t \tilde{\mathcal{B}_t}$ is a (nonsmooth) monotone isotopy - but one whose interfaces is a graph over $\gamma_{d/2}(1)$, so that we may obtain a smooth monotone isotopy $\overline{\mathcal{A}_t \#_t \mathcal{B}_t}$.
$\medskip$ 

One may ask now how this abomination is related to our original $\mathcal{M}$. Now taking the support of the mollifier small enough, we see that $\overline{\mathcal{A}_0 \#_0 \mathcal{B}_0}$ is monotonically isotopic to $\mathcal{M}$ by straightline homotopy since the original ``pinched neck"  (i.e. figure 11) is. We also have that $\overline{\mathcal{A}_1 \#_1 \mathcal{B}_1}$ is two round spheres smoothly glued at the ends of a thin neck.
$\medskip$

So finally we have arrived at step 1 of figure 8 above. The round sphere we note (namely since it is convex) can be written as a graph over the last cross-section of the neck, so using a straightline homotopy again we get a monotone isotopy to step 2 of figure 8. 
 \begin{center}
$\includegraphics[scale = .35]{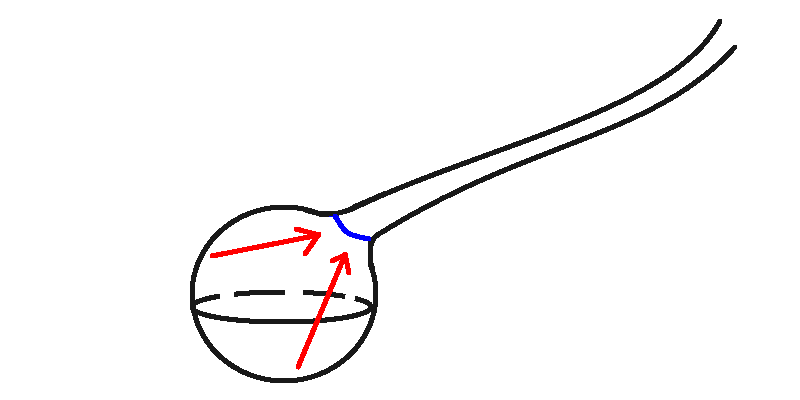}$
\smallskip

$\textbf{Figure 16}$ 
\end{center}

Proceeding as in proposition 5.1 we next find ourselves at step 3; the surface can be written as a graph over a slightly smaller round sphere so using straightline homotopy again we finish. \end{pf}
$\medskip$
\bigskip

\section{Torus to knot theorem}
Now we will prove that every mean convex torus in $\mathbf{R}^3$ is isotopic to an $\epsilon$-thick knot as defined in the introduction. Suppose we have such a torus $\mathcal{M}$. 
$\medskip$

If the flow (with surgeries) of $\mathcal{M}$ has no proper necks in that after the first singular time $T_0$, $M(t) = \emptyset$ (this happens if the neck continuation theorem never stops at caps) then we see that as $t \to T_0$, $M(t)$ converges to a knot (that is, an embedding of $S^1$ in $\R^{n+1}$) in $C^0$ norm so the statement is true (how close depends on how large we took $H_1$ in the surgery procedure \footnote{from \cite{L}, in fact we know as $H_1 \to \infty$ the flow with surgeries converges to the level set flow}).
$\medskip$

\begin{lem} ($\textbf{surgery dichotomy}$) Suppose at surgery time $T_0$ there are standard surgeries; then we can classify performing surgery at a point into one of two cases:
\begin{enumerate}
\item Performing surgery at that point leave a connected manifold diffeomorphic to the sphere.
$\medskip$

\item Performing surgery at that point leaves a disconnected manifold, one diffeomorphic to the torus and the other diffeomorphic to the sphere
\end{enumerate}
\end{lem}
\begin{pf}
The cases in pictures, not including other possible necks of $\mathcal{M}_{T_0}$:
\begin{center}
$\includegraphics[scale = .3]{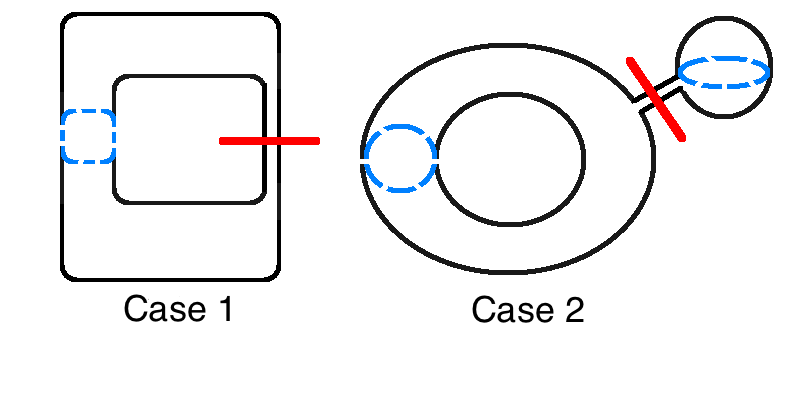}$

$\textbf{Figure 17}$ 
\end{center}
Now to see this dichotomy is true, first suppose the surgery disconnects $\mathcal{M}$, so that $S^{n-1} \times S^1 \cong \mathcal{M}\cong A \# B$ where $A$ and $B$ are both 2-convex. For $n \geq 3$, the seifert van kampen theorem shows up that $\Z = \pi_1(\mathcal{M}) = \pi_1(A)$*$\pi_1(B)$, the free product. Since the free product of two nontrivial groups is never commutative, one of $\pi_1(A)$ or $\pi_1(B)$ is $\{e\}$ and the other is $\Z$, say $\pi_1(A) = \{e\}$. By the classification theorem then $A \cong S^{n}$, and $B = S^{n-1} \times S^1$. 
$\medskip$

If $n = 2$ using the seifert van kampen theorem is more complicated but we can proceed by simpler means anyway. We recall that $0= \chi(\mathcal{M}) = \chi(A \# B) = \chi(A) + \chi(B) - 2$ which implies that $g_\mathcal{M} = g_A + g_B$, so again from our classification theorem one of $A$ or $B$ has to be diffeomorphic to $S^2$ and the other must be diffeomorphic to $S^1 \times S^1$. 
$\medskip$

Now suppose that the surgery leaves $\mathcal{M}$ connected. $\mathcal{M}$ is still 2-convex after the surgery, so it is diffeomoprhic to either $S^n$ or a connect sum $S^{n-1} \times S^1$. Since the surgery keeps $\mathcal{M}$ connected, we see in $\mathcal{M}$ presurgery there is a homotopically nontrivial loop going through the surgery spot. Hence $\mathcal{M}$ post surgery must have fundamental group with at least one less factor $\Z$ ($\pi_1(S^{n-1} \times S^1) = \Z \times \Z$ if $n = 2$ and $\Z$ if $n > 2$). If $n > 2$, then $\mathcal{M}$ post surgery is simply connected so we see from theorem 1.1 (aka corollary 1.3 in \cite{HS2}) it must be diffeomorphic to $S^n$. If $n =2$, none of the surfaces listed have $\pi_1 = \Z$ so $\mathcal{M}$ must be $S^2$ in this case as well. \end{pf}
$\medskip$

Continuing on, if there is a standard surgery as in case 1, then doing surgery here (and only here) leaves a sphere and from theorem 1 (and extending the neck accordingly as in its proof) we have a thin neck attached to a round sphere as below:
\begin{center}
$\includegraphics[scale = .3]{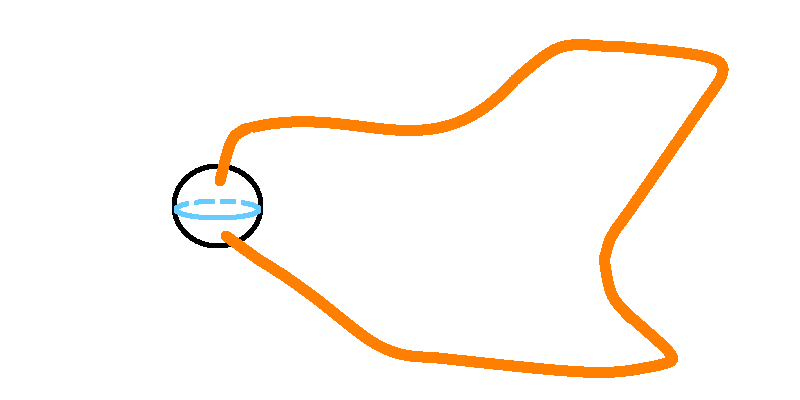}$
\smallskip

$\textbf{Figure 18}$ 
\end{center}
where of course the extended neck (in red) is stylized but can be taken to have uniform diameter as small as we want, by taking $H_1$ large in the surgery definition. Now the next step is to crush the sphere to get to this picture so that we have the result in this case.
$\medskip$

Now suppose though that we only have necks as in case 2 at $T_0$, then we decompose $\mathcal{M}_{T_0}$ as in the induction step in theorem 1's proof (applied to each of the (case 2) necks) and label the pieces as $T$, $B_1, \cdots B_k$ (below $k = 2$) where $T$ as suggested in pictures is diffeomorphic to the torus and the $B_i$ are diffeomorphic to $S^n$:
\begin{center}
$\includegraphics[scale = .3]{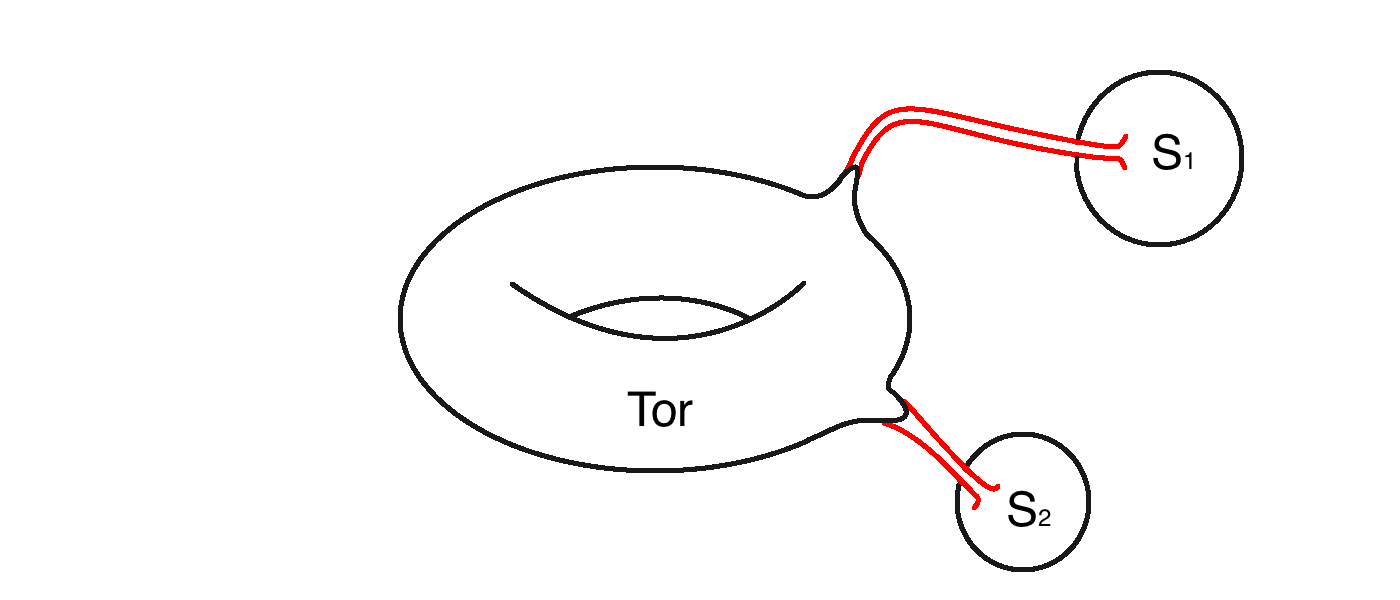}$
\smallskip

$\textbf{Figure 19}$ 
\end{center}
Now by theorem 1, each $B_i$ is isotopic to a round $S^2$. Of course the surgery leaves $T$ mean convex so we may continue to flow it (extending the necks attaching the $B_i$ as in theorem 1) and repeat the process at singular times. As we hit each next singular time (given no case 1 necks) we redefine $T$ to be the ''torus component" (that is, the (would-be) post-surgery component that isn't simply connected - labeled Tor above) left over from the surgeries.
$\medskip$

 If there is ever a neck as in case 1 we get that $T$ can be isotoped to a thick knot. Since the flow with surgeries (untampered with) extinguishes in finite time if there are no case 1 surgeries than $T$ as noted above must flow under its own devices to a thick knot (the ''no surgeries" case). Hence either way we can isotope to the case of several round spheres glued to a thick knot (by a ``hand built neck" as in the sphere connectedness theorem).
 \begin{center}
$\includegraphics[scale = .3]{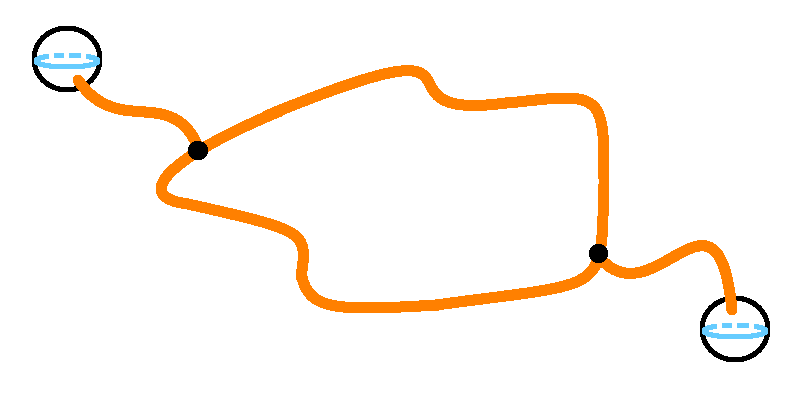}$
\smallskip

$\textbf{Figure 20}$ 
\end{center}
Now as in the end of the proof of theorem 1 (namely, steps 2 -4), we may retract each of the spheres to end up with an $\epsilon$-thick knot, for some knot say $\gamma$. Note that shrinking this further (by straightline homotopy like in the proof of proposition 5.1 above), $\mathcal{M}$ is isotopic to a $\epsilon$ thickening of $\gamma$ for any $\epsilon > 0$, so we are done. 
$\medskip$

\section{Hypersurface to skeleton theorem}
We will essentially use our work and methods from the previous two theorems; this is another proof by induction, but this time we induct on number of tori. From Theorem 1.1, the classification theorem, we know that a 2-convex closed hypersurface $\mathcal{M} \subset \R^{n+1}$ is diffeomorphic to either $S^n$ or a finite connect sum of $S^{n-1} \times S^1$. 
$\medskip$

The cases $\mathcal{M} = S^n$ and $\mathcal{M} = S^{n-1} \times S^1$ are covered respectively in the sphere connectedness theorem (shrink the sphere down even more until its is a sphere of radius $\epsilon$) and the tori to knot theorem. So suppose that $\mathcal{M} \cong (S^{n-1} \times S^1) \# \cdots \# (S^{n-1} \times S^1)$, say $k \geq 2$ direct sums. We proceed by induction on $k$.
$\medskip$

This time there are three cases for standard surgeries, as can be seen similar to above in the torus to knot theorem:
\begin{lem} ($\textbf{surgery trichotomy}$) either performing surgery as a standard surgery point
\begin{enumerate}
\item disconnects $\mathcal{M}$ into $A$ and $B$ where $A$ and $B$ both are connect sums of strictly less copies of $S^{n-1} \times S^1$,
\item  $A \cong S^n$, $B \cong (S^{n-1} \times S^1)^k$
\item leaves $\mathcal{M}$ connected, forcing $\mathcal{M} \cong (S^{n-1} \times S^1)^{k-1}$
\end{enumerate}
\end{lem}
Now in the first and third cases, we may proceed by induction and ``extending the necks" like we have. If at the first surgery time we only have to do standard surgeries of type (2), we continue as in the torus to knot theorem by isotoping the sphere components to round spheres and proceeding on; eventually we will run into standard surgeries of type either 1 or 3, because otherwise like in the finish of the torus to knot theorem $\mathcal{M}^n$ would be diffeomorphic to $S^{n-1} \times S^1$, but we assumed $k \geq 2$. 
$\medskip$

Finally, note that in the tori to knot theorem that ``side" spheres branching off the ``knot" part of the torus were eventually retracted to the knot, so that these sets created are indeed skeletons as described in the introduction.

\section{Finiteness theorem} 

From our $\alpha$ non-collapsed assumption, initial upper bound on $H$, and monotonicity assumption we see there is a uniform lower bound $\delta$ on the diameter of the maximal tubular neighborhood of the skeletons of manifolds in our class $\Sigma(d,C,\alpha)$. Now cover $B^{n+1}_d(0)$ with closed sets $C_i$, disjoint except possibly at their boundaries, with diam$(C_i) \leq \delta/3$ given by intersecting $B_d(0)$ with cubes of side lengths $\ell < \frac{\delta}{3\sqrt{n}}$
$\medskip$ 

With this covering in mind, take $\mathcal{M} \in \Sigma(d, C, \alpha)$ and let $\gamma$ be its corresponding skeleton provided by the hypersurface to skeleton theorem. Since $\gamma$ has at least a $\delta$ thick tubular neighborhood, we may isotope $\gamma$ to another embedded curve (which we will also call $\gamma$) so that $\gamma$ intersects $\partial C_i$ (for any $i$) only through faces, and from the tubular neighborhood assumption we also see that $\gamma_i :=\gamma \cap C_i$ if nonempty is either
\begin{enumerate}
 \item an embedded curve $\gamma_i$ or,
 \item an $\textbf{intersection node}$, where $\gamma$ locally looks like a collection of embedded curves $\gamma_{i_k}$ (with boundary) intersecting at a single point $p_i$.
 \end{enumerate}
  By the tubular neighborhood assumption and that diam$(C_i) < \delta/3 < \delta$ the $\gamma_{i_k}$ must not be tangled and no two $\gamma_{i_k}$ may intersect the same face of $\partial C_i$: 
\begin{center}
$\includegraphics[scale = .28]{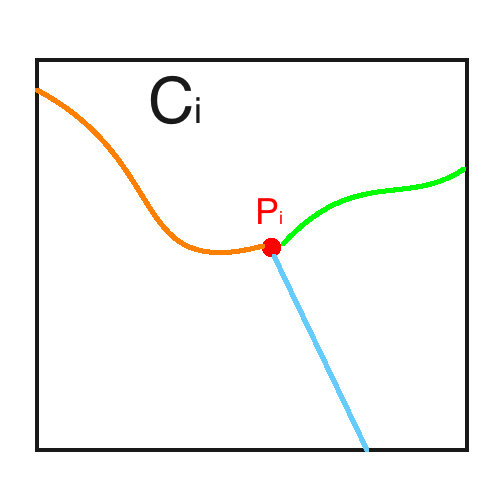}$

$\textbf{Figure 21}$ 
\end{center}
Choose a labeling $i_k \in \{1, \ldots, 2n \}$ of the faces $F_{i_k}$ of $C_i$ so that $\gamma_{i_k}$ intersects $\partial C_i$ on $F_{i_k}$. From our ``tubular assumption" we see we can isotope each of the $\gamma_{i_k}$ so that they intersect $F_{i_k}$ at some distinguished points $D_{i_k} $ in the interior of $F_{i_k}$ which we may pick and so that at every time along the isotopy, $\gamma_{i_k}$ intersects $F_{i_k}$ (we don't detach the segment from the face). Again using the assumption we may isotope them to straight line segments connecting $D_{i_k}$ and $P_i$. Then clearly we can isotope each of the curves so that they intersect at some distinguished point $Q_i$ in the interior of $C_i$ which again we may pick so that we get:
\begin{center}
$\includegraphics[scale = .3]{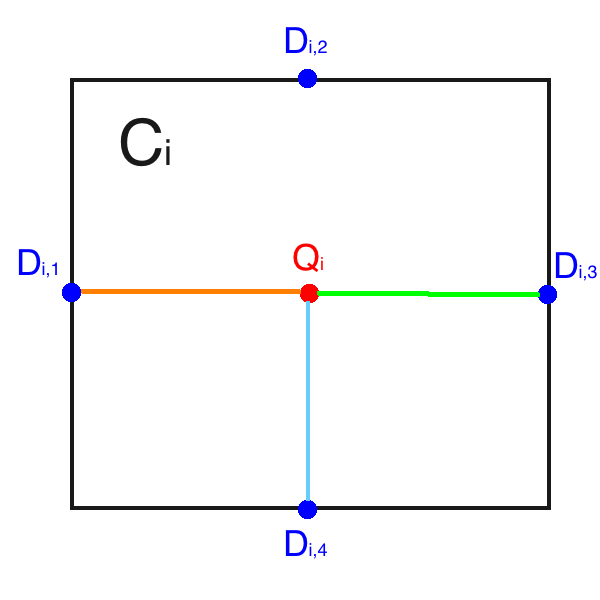}$
\smallskip

$\textbf{Figure 22}$ 
\end{center}

Now suppose that $\gamma_i = C_i \cap \gamma$ isn't an intersection node (or empty). Again choosing a labeling $i_k$ of the faces $F_{i_k}$ of $C_i$, suppose without loss of generality that $\gamma_i$ intersects $\partial C_i$ at $F_{i_1}$ and $F_{i_2}$ (recall skeletons don't have boundary so if $ C_i \cap \gamma$ is nonempty $\gamma$ has to ``leave" $C_i$). Like above then we can isotope $\gamma_i$ so that it intersects $D_{i_1}$ and $D_{i_2}$. Since $C_i$ is contractible, we can also isotope $\gamma_i$ so that it intersects any distinguished point $Q_i$ in the interior of $C_i$ and is in fact the union of two straight lines connecting $D_{i_1}$ and $D_{i_2}$ to $Q_i$. 
$\medskip$

Now for our (closed) cubical cover $\{C_i\}$ of $B_d(0)$ pick any points $Q_i \in \text{int}(C_i)$ and $D_{i_k} \in \text{int}(F_{i_k})$. Moreover make this selection coherent in that if $\partial C_i \cap \partial C_j$ contains a face then their corresponding distinguished points $D_{j}$ coincide (its clearly possible to do this by an iterative process). Then we see from our construction and our coherent choice of $D_{i_k}$ that the isotopies of the $\gamma_i$ described above can be done coordinated to give an isotopy of $\gamma$, so we get something as follows:

\begin{center}
$\includegraphics[scale = .3]{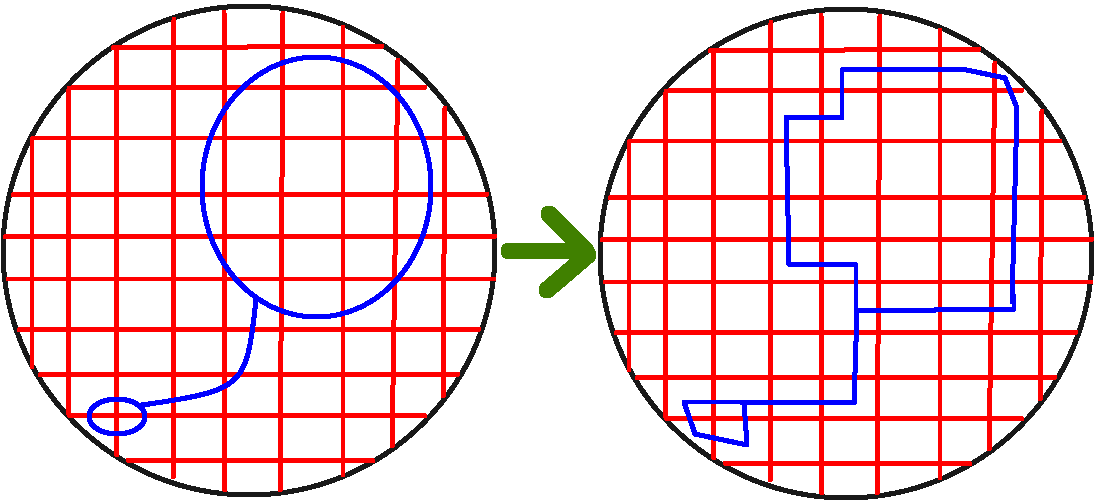}$
\smallskip

$\textbf{Figure 23}$ 
\end{center}
The path on the right we note is only determined by the cubes $C_i$ that $\gamma$ intersects and since there are only finitely many possible paths we have that up to isotopy there are only finitely many skeletons of manifolds in $\Sigma(d,C, \alpha)$ save for one thing, in that we made a choice in isotoping $\gamma$ initially to ensure it intersected all the $C_i$ at (open) faces.  But we see the number of choices is also bounded since the number of cubes $C_i$ is bounded.  Now all that remains to show:
\begin{prop} If $\mathcal{M}_0$ and $\mathcal{M}_1$ have isotopic skeletons, then $\mathcal{M}_0$ is isotopic to $\mathcal{M}_1$
\end{prop}
\begin{pf}
Let $\gamma_0$ be a skeleton of $\mathcal{M}_0$ and similarly let $\gamma_1$ be a skeleton of $\mathcal{M}_1$ (note the skeleton of a hypersurface is not unique, you can perturb it slightly). Let $\gamma_t$ be the (image of) the isotopy of $\gamma_0$ to $\gamma_1$. Note there is a uniform lower bound $\delta$ on the diameter of the (maximal) tubular neighborhood $\gamma_1$ by compactness. 
$\medskip$

Now if $\gamma_0$ is a skeleton of $\mathcal{M}_0$, that means for all $\epsilon > 0$, $\mathcal{M}_0$ is isotopic to an $\epsilon$-thickening $\mathcal{M}_{0,\epsilon}$ of $\gamma_0$. If we take $\epsilon << \delta$ from above, we see the isotopy of $\gamma_0$ to $\gamma_1$ gives rise to an isotopy between $\mathcal{M}_{0, \epsilon}$ and a $\epsilon$-thickening of $\gamma_1$. This, in turn, is (since $\gamma_1$ is a skeleton of $\mathcal{M}_1$) isotopic to $\mathcal{M}_1$, giving the statement. \end{pf}
$\medskip$

Let's try to get a concrete upper bound for this number. First, we see $\delta = \alpha/2C$ works, and (using $\delta = \frac{\alpha}{2C}$, $\ell = \frac{\delta}{6\sqrt{n}} =\frac{\alpha}{12C\sqrt{n}}$ and since $B_d$ is contained in cube of sides length $d$) then we can bound the number of cubes in our cubical cover by
\begin{center}
$\frac{\text{ volume of } B_d }{ \text{volume of each cube}}  \leq \frac{d^n}{ \ell^n} = \frac{(12dC\sqrt{n})^n}{\alpha^n}$

 \end{center}
 For each cube there are $2n$ sides, so since there are at most $2n$ segements of a skeleton leaving a cube in the cover (with our chosen $\delta$) from the discussion above. Hence we can bound the number of skeletons up to isotopy by $2^{2n} \frac{( 12dC\sqrt{n})^n}{\alpha^n}$
$\medskip$

\section{Concluding Remarks}

The upper bound provided in the finiteness theorem is quite crude. For one, we count the number of cubes in the cubical cover just by counting how many cubes it takes to cover a cube that contains of side length $d$ which is clearly overcounting. Secondly, as discussed in the proof this is some room for ambiguity in that several ``standard" skeletons may correspond to the same class up to isotopy. Our counting also includes nonconnected skeletons, which we aren't really interested in. It seems it would be interesting to know though, letting $d \to \infty$, $\alpha \to 0$, or $C \to \infty$ what the correct growth rate is asymptotically and perhaps this is achieved with a sloppy counting scheme of our type. 
$\medskip$

An upper bound on genus is also implied by the finiteness theorem, but this seems even harder to calculate with any degree of accuracy using a simple counting method. It might be interesting though, using techniques from random graph theory, to use the reduction of hpyersurfaces to skeletons to make a statement on what the ``average" genus of a hypersurface in $\Sigma(d,C,\alpha)$ is. 
$\medskip$

The procedures described above seem to not depend continuously on initial hypersurface $\mathcal{M}$. To see why, consider the an isotopy of a sphere to a dumbbell by ``pulling" the sphere apart like in the below figure (parameterized by $s$):
\begin{center}
$\includegraphics[scale = .3]{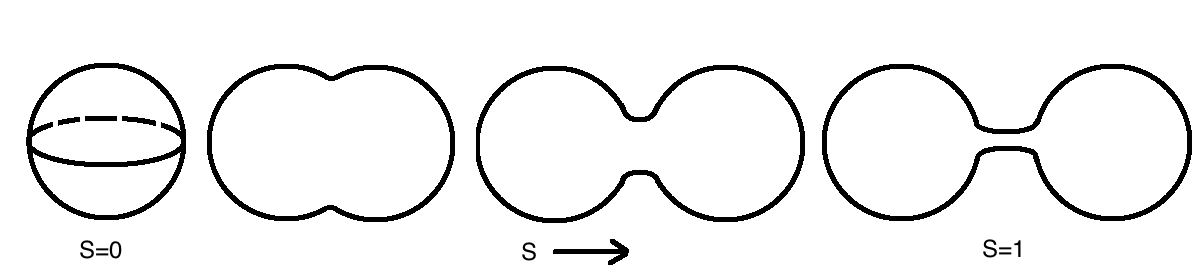}$
\smallskip
\end{center}
The dumbbell ($s=1$) is designed to develop a neckpinch singularity where the bells are low curvature at the singular time. Of course when $s=0$ the round sphere just shrinks in on itsself and under the algorithm described in the path to sphere theorem will flow by the mean curvature flow untill its curvature is high and will be straightline homotoped to a small round sphere (of course in this case, the theorem is trivial!). We see then that there is an intermediate time $s =s_0$ when for $s < s_0$, the surface will eventually become high curvature all at once and for $s \geq s_0$ the hypersurface will develop a neckpinch singularity, with a neck bordered by low curvature regions (to make sure the high curvature region is bordered by $\textbf{two}$ low curavature regions, arrange things so they are reflection symmetric). In that case, where the surgery would have been done the algorithm described ``straightens" the neck as in figure 10. In the other case for $s \leq s_0$ the straightline isotopy involved is nonconstant on the whole hypersurface for all $t \in [0,1]$. One sees then that the algorithm is discontinuous at $s = s_0$ (at $s = s_0$ we should have a ``degenerate neckpinch" ) . This discussion hopefully also illustrates the subtlety of the issue. 
$\medskip$

To extend these ideas to a wider class of hypersurfaces of course seems to boil down to extending the class of submanifolds on which one can do surgery. There are three main ways one imagines to move forward; perhaps most naturally is to try to the relax the curvature conditions,  perform surgery where the ambient manifold isn't Euclidean space, and also extend the surgery to higher codimension submanifolds:
$\medskip$

Relaxing the curvature conditions past mean convexity seems to lead to complications for this paper's argument, because the monotonicity of the flow will no longer hold; one could imagine though by ``moving the necks out of the way" analogous to some arguments found in \cite{BHH} along the flow it could be possible to still pass the argument through. Much more seriously is that the zoo of possible singularities (and hence high curvature regions) arising under the flow broadens considerably in the general case and its not clear one can control their topology; for example compact genus one singularities are possible (see \cite{Ang}) and there is strong numerical evidence even compact higher genus examples exist (see \cite{Cho} and also \cite{Ket} where (possibly noncompact) self shrinkers with the same symmetries are constructed by min-max methods). Even in the mean convex case, where singularities are either round $S^n$ or $S^{n-k} \times \R^k$, it is unclear how to handle the case of a ``sheet" of singularities (when $k > 1$). 
$\medskip$

There has been very interesting work recently on extending the surgery methods to different ambient manifolds $N$. From the previous discussion one imagines its most natural to start with the class of two convex hypersurfaces $M$ in $N$, but this condition isn't always preserved by the mean curvature flow so one is naturally drawn to consider other flows. Taking this approach, in \cite{BH1} Brendle and Huisken reprove all the estimates necessary to do surgery for hypersurfaces moving with normal velocity $G_\kappa = \left(\sum\limits_{i < j} \frac{1}{\lambda_i + \lambda_j - 2\kappa} \right)^{-1}$, where $\overline{R}_{1313} + \overline{R}_{2323} \geq - 2\kappa^2$ at each point of $M_t$ - of course this flow preserves 2-convexity ($G_k$ is a particular flow that satisfies certain properties to ensure this, but there are conceivably others -see remark 1.3 in \cite{BH1}). One imagines most if not all of the statements in this paper are true in ambient manifolds using this flow as well. 
$\medskip$

Concerning generalizing past the codimension one case of course the same issues discussed above (namely understanding of singularities) are generally true when the codimension is increased; in higher codimension the mean curvature flow is in fact even less well understood than in the general codimension one case because here the mean curvature is a vector, not just a function. At any rate, a natural curvature condition to impose (although there may be others -see \cite{Bake}) is the Lagrangian condition; it might still be possible to carry out such a surgery program for Lagrangian mean curvature flow due to the special geometry of Lagrangians - for a discussion of this possibility see \cite{JLT}. Looking further down the road, recall Nash's isometric embedding theorem, that any manifold can be isometrically embedded into $\mathbf{R}^N$ for some $N$ but it might be of high codimension. So, if we could extend the surgery to manifolds of arbitrarily high codimension with no assumptions on the curvature, and if we could still run the argument above, then this seems to preclude the existence of exotic spheres! So there seems to be fundamental limitations to fulfilling ones hopes in this direction. 
$\medskip$

\end{document}